\documentclass[11pt]{article}
\usepackage{latexsym}
\usepackage{amsmath}
\usepackage{color}
\usepackage{graphicx}

\newtheorem{thm}{Theorem}[section]
\newtheorem{lem}[thm]{Lemma}
\newtheorem{cor}[thm]{Corollary}

\newtheorem{rem}[thm]{Remark}
\newcommand{\AS}{{\rm a.s.}}
\newcommand{\Exp}{{\rm E}}
\newcommand{\Var}{{\rm Var}}
\newcommand{\convweak}{{\Rightarrow}}
\newcommand{\convprob}{\stackrel{P}{\rightarrow}}

\newcommand{\eqdist}{\stackrel{\rm d}{=}}
\newcommand{\bbR}{{\rm I\hspace{-0.8mm}R}}
\newcommand{\matR}{{\bbR}}
\newcommand{\koniec}{\newline\vspace{3mm}\hfill $\Box$}

\title{Strong approximations for long memory sequences based
partial sums, counting and their Vervaat processes
\protect}
\author{
Endre Cs\'{a}ki\thanks{Research supported by the Hungarian National
Foundation for Scientific Research No. K108615}\\
A. R\'enyi Institute of Mathematics, Hung. Acad. Sci. \\
Budapest, P.O.B. 127, H-1364 Hungary\\
e-mail: csaki.endre@renyi.mta.hu
\and
Mikl\'os Cs\"org\H o\thanks{Research supported by an NSERC Canada
Discovery Grant at Carleton University, Ottawa}\\
School of Mathematics and Statistics, Carleton University\\
Ottawa, Ontario, Canada K1S 5B6\\
e-mail: mcsorgo@math.carleton.ca
\and
Rafa{\l} Kulik\thanks{Research supported by an NSERC Canada Discovery
Grant at the University of Ottawa, Ottawa}\\
Department of Mathematics and Statistics, University of Ottawa\\
Ottawa, Ontario, Canada K1N 6N5\\
e-mail: rkulik@uottawa.ca
}

\begin{document}

\maketitle
\begin{quotation}
\noindent {\it Abstract:}
We study the asymptotic behaviour of partial sums of long range dependent
random variables and that of their counting process, together with an
appropriately normalized integral process of the sum of these two
processes, the so-called Vervaat process. The first two of these processes
are approximated by an appropriately constructed fractional Brownian
motion, while the Vervaat process in turn is approximated by the square of
the same fractional Brownian motion.
\par

\vspace{9pt}
\par

\end{quotation}\par


\fontsize{10.95}{14pt plus.8pt minus .6pt}\selectfont

\setcounter{equation}{0}

\section{Introduction}
Let $\{\tilde\eta_j,j\ge 0\}$ be a stationary long range dependent (LRD)
centered Gaussian sequence with $\Exp(\tilde\eta_0^2)=1$ and covariance
function of the form
\begin{equation}
\rho_k:=\Exp(\tilde\eta_0\tilde\eta_k)=k^{-\alpha}L(k),\quad
k=1,2,\ldots,
\label{cov}
\end{equation}
where $\alpha\in (0,1)$, and $L(\cdot)$ is a slowly varying function at
infinity. Let $G(\cdot)$ be a real valued Borel measurable function with
$\Exp(G(\tilde\eta_0))=\mu$ and $\Exp(G^2(\tilde\eta_0))<\infty$.
It may be expanded as
$$
G(\tilde\eta_0)-\mu=\sum_{q=m}^{\infty}\frac{J_q}{q!}H_q(\tilde\eta_0),
$$
where convergence is in $L^2$,
$$
H_q(x):=(-1)^qe^{x^2/2}\frac{d^q}{dx^q}e^{-x^2/2}, \qquad
q=1,2,\ldots,
$$
are Hermite polynomials, and $J_q:=\Exp
(G(\tilde\eta_0)H_q(\tilde\eta_0))$. The
index $m$ in the above expansion is defined as
$$
m:=\min\{q\geq 1:J_q\not=0\}.
$$
We then say that the Hermite rank of $G(\cdot)$ is $m$.

Consider the subordinated sequence $\{Y_j=G(\tilde\eta_j),\, j\geq 0\}$.
Then $\Exp(Y_j)=\mu$, $j=0,1,2,\ldots$ Given (\ref{cov}) with
$\alpha\in (0,1)$, assume that the Hermite
rank $m$ of $G(\cdot)$ is such that $0<\alpha<1/m$. Then, as $n\to\infty$,
with $L(\cdot)$ as in (\ref{cov}), we have
(cf. Lemma 3.1 and Theorem 3.1 of Taqqu \cite{Taqqu1975}, or
\cite{Taqqu1977}),
\begin{equation}
\sigma_{n,m}^2:=\Var \left(\sum_{j=1}^n(Y_j-\mu)\right)\sim
\frac{J_m^2}{m!}\frac{2}{(1-\alpha m)(2-\alpha m)}n^{2-m\alpha}L^m(n),
\label{sigma}
\end{equation}
where the symbol $\sim$ stands for the indicated terms being
asymptotically equal to each other.

We note that an LRD Gaussian sequence as in (1) can be viewed in terms
of the linear (moving average) process
\begin{equation}
\eta_j=\sum_{k=0}^\infty \psi_k \xi_{j-k},\qquad j=0,1,2,\ldots,
\label{linear}
\end{equation}
where $\{\xi_k,\, -\infty<k<\infty\}$ is a double sequence of independent
standard normal random variables, and the sequence of weights
$\{\psi_k, \, k=0,1,2,\ldots\}$ is square summable. Then
$\Exp(\eta_0)=0$, $\Exp(\eta_0^2)=\sum_{k=0}^\infty\psi_k^2=:
\sigma^2$ and, on putting $\tilde\eta_j=\eta_j/\sigma$, $\{\tilde\eta_j,\,
j=0,1,2,\ldots\}$ is a stationary Gaussian sequence with
$\Exp(\tilde\eta_0)=0$ and $\Exp(\tilde\eta_0^2)=1$. If
$\psi_k\sim k^{-(1+\alpha)/2}\ell(k)$ with a slowly varying function,
$\ell(k)$, at infinity, then (cf. Wang {\it et al.} \cite{WangLinGulati2003})
$\Exp(\eta_j\eta_{j+n})\sim b_\alpha n^{-\alpha}\ell^2(n)$, i.e.,
we have (\ref{cov}) with $L(n)\sim b_\alpha\ell^2(n)/\sigma^2$, where the
constant $b_\alpha$ is defined by
$$
b_\alpha=\int_0^\infty x^{-(1+\alpha)/2}(1+x)^{-(1+\alpha)/2}\, dx.
$$

For more information on long memory Gaussian sequences and their
functionals we refer to \cite{Taqqu1975} and \cite{Taqqu2003}.

In this exposition we study the asymptotic behaviour of partial sums of
long range dependent random variables and that of their counting process,
together with an appropriately normalized integral process of the sum of
these two processes, the so-called Vervaat process. All this will amount
to an extension of the recent work in \cite{CsakiCsorgoKulik} and
\cite{CsakiCsorgoRychlikSteinebach2007} to LRD sequences.

\subsection{Bahadur-Kiefer processes}
For initiating our discussion, let $U_j,\, j=1,2,\ldots$, be
independent identically distributed uniform random variables on $[0,1]$,
and denote by $F_n$ the right continuous empirical distribution function
of the first $n\ge 1$ of these random variables. Define the so-called
uniform Bahadur-Kiefer process $\{R_n(t), 0\le t\le 1\}$ as
\begin{equation}\label{eq:BahKiefer}
R_n(t):=\alpha_n(t)+\beta_n(t),
\end{equation}
where $\alpha_n(t):=n^{1/2}(F_n(t)-t)$,
$\beta_n(t):=n^{1/2}(F_n^{-1}(t)-t)$, and $F_n^{-1}(\cdot)$ is the
empirical quantile function (left continuous inverse of $F_n(\cdot)$). The
process $R_n(\cdot)$ was introduced and first studied by Bahadur
\cite{Bahadur1966} and, consequently, by Kiefer in \cite{Kiefer1967}
and \cite{Kiefer1970}. In particular, it can be concluded from the results
of Kiefer in the just cited papers that the Bahadur-Kiefer process
$R_n(\cdot)$ \textit{cannot converge weakly to any non-degenerate
random element of the space $D[0,1]$}. This can also be deduced via
the following result of Vervaat in  \cite{Vervaat1972a} and
 \cite{Vervaat1972b}: as $n\to\infty$,
\begin{equation}\label{eq:Vervaat}
V_n(t):=2n^{1/2}\int_0^tR_n(u)\, du\, \convweak\, W_0^2(t),\qquad 0\le t\le 1,
\end{equation}
where $W_0(t)$ is a standard Brownian bridge, and $\convweak$ stands
for weak convergence in $C[0,1]$, equipped with the uniform norm. Since
$W_0^2$ is not differentiable, we conclude again that $R_n$, as the
derivative of $V_n$ cannot converge weakly in $D[0,1]$ with any normalization.
The integrated Bahadur-Kiefer process in  (\ref{eq:Vervaat}) is usually
referred to as the uniform Vervaat process.

The just mentioned works of Vervaat constitute important general
contributions to limit theorems for processes with positive drift and
their inverses which we now indicate briefly. Accordingly, let ${\cal Z}$
be a non-negative stochastic process on $[0,\infty)$ such that almost all
realizations ${\cal Z}(\cdot,\omega):[0,\infty)\to [0,\infty)$ are
non-decreasing unbounded functions, and ${\cal Z}^{-1}$ be the generalized
inverse of ${\cal Z}(\cdot,\omega)$, i.e., ${\cal Z}^{-1}(t,\omega):=
\inf\{u:\, {\cal Z}(u,\omega)>t\}$. Let $D_0[0,\infty)$ be the subset of
non-decreasing, non-negative unbounded functions of $D[0,\infty)$, the set
of real-valued functions on $[0,\infty)$ which are right-continuous and
have finite left-hand limits at every point $t\in (0,\infty)$, and
let $C[0,\infty)$ be the subset of continuous functions of
$D[0,\infty)$. For further use we state here Theorem 1 of Vervaat
\cite{Vervaat1972a} for our convenience as follows (cf. also Theorem
3.2.3 of Vervaat \cite{Vervaat1972b}). Here and throughout $\convweak$
indicates weak convergence in an appropriate context, while $\convprob$
designates convergence in probability.

\begin{thm}\label{thm:Vervaat}
Let ${\cal Z}_1,{\cal Z}_2,\ldots,$ be random elements in $D_0[0,\infty)$,
$\tilde {\cal Z}$ a random element in $C[0,\infty)$, and
$\zeta_1,\zeta_2,\ldots, $ be positive random variables such that
$\zeta_n\convprob 0$ as $n\to\infty$. Then, as $n\to\infty$, the following
two weak convergence statements are equivalent in $D[0,\infty)$ (endowed
with the uniform topology on compact sets):
$$
\frac{{\cal Z}_n-I}{\zeta_n}\, \convweak\, \tilde {\cal Z},
$$
$$
\frac{{\cal Z}_n^{-1}-I}{\zeta_n}\, \convweak\, -\tilde {\cal Z},
$$
where $I$ denotes the identity map on $[0,\infty)$. Moreover, if any
of the above statements holds, then
\begin{equation}
V(\cdot\, ;{\cal Z}_n)\, \convweak\, \frac{1}{2}\widetilde {\cal Z}^2
\label{v1}
\end{equation}
in $C[0,\infty)$ (endowed with the uniform topology on compact
sets), where
$$
V(t;{\cal Z}_n)=\frac{1}{\zeta_n^2}
\int_0^t({\cal Z}_n(u)+{\cal Z}_n^{-1}(u)-2u)\, du,
\qquad 0\le t<\infty.
$$
\end{thm}

Clearly, the statement of (\ref{eq:Vervaat}) is implied by that of
(\ref{v1}).

On the other hand, if we consider a subordinated sequence of random
variables $Y_j=G(\tilde\eta_j)$, $j=0,1,2,\ldots,$ with marginal
distribution function $F(x)=P(Y_0\leq x)$, $x\in{\matR}$, where
$\{\tilde\eta_j,\, j\geq 0\}$ is a stationary LRD centered Gaussian
sequence with $E(\tilde\eta_0^2)=1$ and covariance as in (\ref{cov}), then
their appropriately scaled Bahadur-Kiefer process {\it does converge weakly}.
To be more specific, on expanding $1\{Y_0\le x\}-F(x)$ in terms of Hermite
polynomials for any fixed $x\in {\matR}$ as in Dehling and Taqqu
\cite{DehlingTaqqu1989b}, we have
$$
1\{Y_0\le x\}-F(x)=\sum_{q=m_x}^\infty \frac{c_q(x)}{q!}H_q(\tilde\eta_0),
\quad x\in {\matR},
$$
where $c_q(x)=\Exp (1\{G(\tilde\eta_0)\le x\}-F(x))H_q(\tilde\eta_0)$,
and $m_x$ for any $x\in {\matR}$ is the index of the first nonzero
coefficient in this expansion, the so-called Hermite rank of
$1\{Y_0\le x\}-F(x)$. Then, the Hermite rank of the class
$(1\{Y_0\le x\}-F(x),\, x\in {\matR})$ is defined as
\begin{equation}
m:=\min\{m_x:c_{m_x}(x)\not=0\quad\mbox{\rm for some } x\in {\matR}\}.
\label{m}
\end{equation}

On assuming now that $F(\cdot)$ is continuous, in terms of the first
$n\geq 1$ of the uniformly distributed random variables
$U_j=F(Y_j)=F(G(\tilde\eta_j))$, $j=1,2,\ldots$, keeping the notation we
already used in the i.i.d. case in (\ref{eq:BahKiefer}) for the LRD case
as well, we redefine $\alpha_n(\cdot)$ and $\beta_n(\cdot)$, respectively, as
$$
\alpha_n(t):=\frac{n}{d_{n,m}}(F_n(t)-t),\qquad t\in [0,1],
$$
$$
\beta_n(t):=\frac{n}{d_{n,m}}(F_n^{-1}(t)-t),\qquad t\in [0,1],
$$
where $d_{n,m}^2=n^{2-m\alpha}L^m(n)$.

In this context, from Theorem 1.1 of Dehling and Taqqu
\cite{DehlingTaqqu1989b}, we obtain that as long as
$0<\alpha<1/m$, as $n\to\infty$,
\begin{equation}\label{eq:Dehling-Taqqu}
\alpha_n(t)\, \convweak\,
\sqrt{\frac{2}{(2-m\alpha)(1-m\alpha)}}c_m(F^{-1}(t))X_m,
\end{equation}
in $D[0,1]$, equipped with the sup-norm, where $F^{-1}(\cdot)$ is the
quantile function (inverse) of $F(\cdot)$, $X_m\eqdist X_m(1)$, where
$\{X_m(s),\, 0\leq s \leq 1\}$ is $1/m!$ times a Hermite process of rank
$m$, given for each $s$ as a multiple Wiener-It\^o-Dobrushin
integral (see, e.g., \cite{DehlingTaqqu1989b}, \cite{Taqqu2003}).

As a consequence of (\ref{eq:Dehling-Taqqu}) and Theorem
\ref{thm:Vervaat}, we immediately conclude also that, as $n\to\infty$,
$$
\beta_n(t)\, \convweak\,
-\sqrt{\frac{2}{(2-m\alpha)(1-m\alpha)}}c_m(F^{-1}(t))X_m,
$$
in $D[0,1]$, equipped with the sup-norm.

Furthermore, on using the same definition of the Bahadur-Kiefer process
as in (\ref{eq:BahKiefer}) in our present LRD context as well, via
(\ref{v1}) we arrive at
\begin{equation}\label{eq:Vervaat-LRD}
\frac{n}{d_{n,m}}\int_0^tR_n(u)\, du\, \convweak\,
\frac{1}{(2-m\alpha)(1-m\alpha)}c_m^2(F^{-1}(t))X_m^2,
\end{equation}
in $C[0,1]$, endowed with the uniform topology.

Now, unlike in the i.i.d. case (cf. (\ref{eq:Vervaat})), we can formally
differentiate the processes in (\ref{eq:Vervaat-LRD}), and thus obtain the
following corollary.
\begin{cor}
Let $0<\alpha<1/m$, where $m$ is as in $(\ref{m})$, and let $X_m$ be the
random variable defined in $(\ref{eq:Dehling-Taqqu})$. Assume that the
marginal distribution function $F$ has a positive density $f=F'$, with
respect to Lebesgue measure, on $(a,b)$, where $a=\sup\{x:F(x)=0\}$,
$b=\inf\{x:F(x)=1\}$, $-\infty\le a<b\le\infty$. Assume that the function
$c_m(\cdot)$ is also differentiable over this interval $(a,b)$. Then, as
$n\to\infty$,
\begin{equation}\label{eq:BahKiefer-LRD}
\frac{n}{d_{n,m}}R_n(t)\, \convweak\,
\frac{2}{(2-m\alpha)(1-m\alpha)}
\frac{c_m(F^{-1}(t))c_m'(F^{-1}(t))}{f(F^{-1}(t))}X_m^2,
\end{equation}
in $D[0,1]$, endowed with the sup-norm, provided the deterministic function
$$\frac{c_m(F^{-1}(t))c_m'(F^{-1}(t))}{f(F^{-1}(t))}$$
is finite for $t\in [0,1]$.
\end{cor}

Thus, unlike in the i.i.d. case, the LRD-based
Bahadur-Kiefer process does converge weakly. Dealing with the sequential
version of (\ref{eq:BahKiefer-LRD}), this was first observed in
\cite{CsorgoSzyszkowiczWang2006}; see also \cite{CSW2006cor},
\cite{CsorgoKulik2008a}, \cite{CsorgoKulik2008b}. In these papers the
method of establishing (\ref{eq:BahKiefer-LRD}) was in fact
completely different. Namely, strong approximations were used instead of
making use of Vervaat's approach in hand.
This, however, led to stronger assumptions on $F$. Nevertheless, the
method of strong approximations as in the above mentioned works is needed
to deal with weak and strong laws of the general (non-uniform) and
possibly sequential Bahadur-Kiefer and Vervaat processes that are based
on subordinated sequences or linear processes.

\subsection{Partial Sums and Counting Processes}
We are to see in this exposition that, unlike Bahadur-Kiefer processes,
i.i.d. and LRD based sums of partial sums and their counting processes
behave similarly in that {\sl both} cannot converge weakly.

Specifically now, assume that $\Exp (Y_0)=\mu$, put $S_n:=Y_1+\cdots+Y_n$,
$n\geq 1$, and define, as in
\cite{CsakiCsorgoRychlikSteinebach2007},
\begin{equation}\label{eq:defin-Sn}
S_n(t):=(n\mu)^{-1}S_{[nt]},\qquad N_n(t):=N(n\mu t)/n,
\end{equation}
where $N(t):=\min\{n\ge 1:S_n>t\}$.

Assume first that the random variables $Y_j$ are i.i.d., and, via
$S_n(\cdot)$ and $N_n(\cdot)$ as in (\ref{eq:defin-Sn}) define now the
following analogue of the Bahadur-Kiefer process $R_n$ as in
(\ref{eq:BahKiefer}) as follows:
$$
R_n^*(s):=n^{1/2}(S_n(s)-s+N_n(s)-s).
$$
Then, in view of the left-hand side of the statement in
(\ref{eq:Vervaat}), we define the following Vervaat-type process:
$$
\tilde V_n(t):=n\int_0^t\left\{(S_n(s)-s)+(N_n(s)-s)\right\}\, ds
=n^{1/2}\int_0^tR_n^*(s)\, ds.
$$

In \cite{CsakiCsorgoRychlikSteinebach2007} strong approximations were
established for $V_n(\cdot)$ in terms of the square of a standard Brownian
motion under the conditions $P(Y_0>0)=1$ and $\Exp (Y_0^4)<\infty$.

In \cite{CsakiCsorgoKulik} the condition $\Exp(Y_0^4)<\infty$ was replaced
by a weaker moment condition and the results of
\cite{CsakiCsorgoRychlikSteinebach2007} were extended to weakly
dependent cases.

In view of Donsker's theorem for partial sums of i.i.d. random
variables, it follows from  Theorem 1.1 that, if $P(Y_0\geq 0)=1$ and
$0<\Exp(Y_0^2)<\infty$, the i.i.d. based Vervaat-type process
$\{\tilde V_n(t),t\in [0,1]\}$ converges weakly in $C[0,1]$ (equipped with
the uniform topology) to $\{W^2(t),t\in [0,1]\}$, where the latter process is
the square of a standard Brownian motion. Consequently, similarly to
(\ref{eq:Vervaat}), the ``derivative" of $\tilde V_n(t)$, i.e., the i.i.d.
based integrand process $\{R_n^*(t),t\in [0,1]\}$ as above cannot converge
weakly.

It will be shown in this paper that in case of LRD sequences, for
$R^*_n(\cdot)$ we have a similar lack of weak convergence.  Namely, a
strong approximation of the appropriately scaled LRD sequences based
$\{\tilde V_n(t),t\in [0,1]\}$ yields a limiting process in
terms of the square of a fractional Brownian motion. Consequently, its LRD
based integrand process $\{R_n^*(t),t\in[0,1]\}$ cannot converge weakly.
This is in contrast to having the weak convergence in the LRD case of
the appropriately scaled Bahadur-Kiefer process $R_n(\cdot)$ as in
(\ref{eq:BahKiefer-LRD}).

In order to conclude our strong approximation for an LRD sequences based
Vervaat process, we first establish some new results on the almost sure
behaviour of similarly based partial sums and their counting process.

The paper is organized as follows. In Section \ref{sec:statement}, we
introduce our basic assumptions and notations for the sake of stating
our strong approximation results for LRD based partial sums, their
counting process and the therein defined Vervaat process. In Section
\ref{sec:consequences}, we list weak convergence laws and laws of the
iterated logarithm (LIL) for these LRD based processes via well-known laws
for appropriate fractional Brownian motions. The proofs are given in
Section \ref{sec:proofs}. In particular, in Section \ref{sec:proofmain1},
we first describe the construction of the fractional Brownian motion that
will be used throughout later on when proving our strong approximation
results for our LRD sequences based processes in hand. The latter
approximations will be achieved by an appropriate reduction
principle and strong approximation for partial sums of LRD based
subordinated sequences that, in turn, conclude Section
\ref{sec:proofmain1}. Following this, in Sections \ref{sec:counting} and
\ref{sec:3-3} respectively, we establish our strong approximations for the
counting and Vervaat processes in hand.

\section{Results}\label{sec:results}

\subsection{Statements of results}\label{sec:statement}

Put $\tilde\eta_j=\eta_j/\sigma$, $j=0,1,2,\ldots$, where $\eta_j$ is as
in (\ref{linear}) with $\sigma^2=\Exp(\eta_0^2)=\sum_{k=0}^\infty
\psi_k^2$. 
We shall assume for simplicity that the slowly varying function in (\ref{cov}) is constant, equal to one. 
We are now concerned with the subordinated sequence
$Y_j=G(\tilde\eta_j)$, $j=0,1,2,\ldots$. We assume throughout that
$J_1:=\Exp(G(\tilde\eta_0)\tilde\eta_0)\not=0$, i.e., the Hermite rank
of $G(\cdot)$ is assumed to be 1.

We will make use of, or a part of, the following assumptions:
\begin{itemize}
\item[{\rm (i)}] $\Exp Y_0=\mu>0$;
\item[{\rm (ii)}] $\Exp(Y_0^{2})<\infty$;
\item[{\rm (iii)}] $P(Y_0\ge 0)=1$.
\end{itemize}

In terms of our subordinated sequence $\{Y_j,\, j=0,1,2,\ldots\}$, with
$t\geq 0$, we define
\begin{eqnarray}
S(t)&:=&\sum_{i=1}^{[t]}Y_i,\label{st}\\
N(t)&:=&\inf\{s\geq 1:\, S(s)>t\},\label{nt}\\
Q(t)&:=& S(t)+\mu N(\mu t)-2\mu t,\label{qt}\\
Z(t)&:=& \mu\int_0^t Q(s)\, ds.
\label{zt}
\end{eqnarray}

In terms of these definitions we have the following analogue of the
Bahadur-Kiefer process
\begin{equation}
R_n^*(s):=\frac{S(ns)+\mu N(\mu ns)-2\mu ns}{ n^{1-\alpha/2}}=
n^{\alpha/2}\frac{Q(ns)}{n}.
\label{BK}
\end{equation}
Via $Z(\cdot)$, the Vervaat process in this context is defined by
\begin{equation}
V_n(t):=n^{\alpha/2}\int_0^t R_n^*(u)\, du=n^{\alpha/2}\int_0^t
n^{\alpha/2}\frac{Q(nu)}{n}\, du= \frac{Z(nt)}{\mu n^{2-\alpha}}.
\label{Ve}
\end{equation}

For $1/2<H<1$ let $\{W_H(t),\, t\geq 0\}$  be a fractional Brownian
motion, i.e. a mean-zero stationary Gaussian process with covariance
\begin{equation}
\Exp W_H(s)W_H(t)=\frac{1}{2}(s^{2H}+t^{2H}-|s-t|^{2H}).
\label{covar}
\end{equation}

\begin{thm}\label{thm:main1}
Let $\eta_j$ be defined by {\rm (\ref{linear})} with $\psi_k\sim
k^{-(1+\alpha)/2}$, $0<\alpha<1$, and put $\tilde\eta_j=\eta_j/\sigma$
with $\sigma^2:=\Exp(\eta_0^2)=\sum_{k=0}^\infty \psi_k^2$. Let $G(\cdot)$
be a function whose Hermite rank is $1$. Furthermore, let $\{S(t),\, t\geq 0\}$ 
be as in {\rm (\ref{st})} and assume condition {\rm (ii)}. Then, on
an appropriate probability space for the sequence
$\{Y_j=G(\tilde\eta_j),\, j=0,1,\ldots\}$, one can construct a fractional
Brownian motion $W_{1-\alpha/2}(\cdot)$ such that, as $T\to\infty$, we have
\begin{equation}\label{eq:3/1}
\sup_{0\le t\le T}\left|S(t)-\mu t-
\frac{J_1 \kappa_{\alpha}}{\sigma}W_{1-\alpha/2}(t)\right|=
o_{\AS}(T^{\gamma/2+\delta}),
\end{equation}
where
\begin{equation}
\kappa_{\alpha}^2=
2\frac{\int_0^{\infty}x^{-(\alpha+1)/2}(1+x)^{-(\alpha+1)/2}\,
dx}{(1-\alpha)(2-\alpha)},
\label{kalpha}
\end{equation}
$\gamma=2-2\alpha$ for $\alpha<1/2$,
$\gamma=1$ for $\alpha\geq 1/2$ and $\delta>0$ is arbitrary.
\end{thm}

\begin{thm}\label{thm:main2}
Assume the conditions of {\rm Theorem \ref{thm:main1}} and condition {\rm
(i)}. Then, on the probability space of {\rm Theorem \ref{thm:main1}}
for the sequence $\{Y_j=G(\tilde\eta_j),\, j=0,1,\ldots\}$, together with
the fractional Brownian motion $W_{1-\alpha/2}(\cdot)$ as in {\rm
(\ref{eq:3/1})}, for $N(t)$ as in {\rm (\ref{nt})}, as
$T\to\infty$, we have
\begin{equation}
\sup_{0\leq t\leq T}\left |\mu N(\mu t)-\mu t+\frac{J_1
\kappa_\alpha}{\sigma}
W_{1-\alpha/2}(t)\right
|=o_{\AS}(T^{\gamma/2+\delta}+T^{(1-\alpha/2)^2+\delta})
\label{eq:main2}
\end{equation}
with $\gamma$ as in {\rm Theorem \ref{thm:main1}}, and arbitrary
$\delta>0$.
\end{thm}

\begin{thm}\label{thm:main3}
Assume the conditions of {\rm Theorem \ref{thm:main1}} and condition {\rm
(iii)}. Then, on the probability space of {\rm Theorem \ref{thm:main1}}
for the sequence $\{Y_j=G(\tilde\eta_j),\, j=0,1,\ldots\}$, together with
the fractional Brownian motion $W_{1-\alpha/2}(\cdot)$ as in
{\rm (\ref{eq:3/1})} and {\rm (\ref{eq:main2})}, for $Z(t)$ as in
{\rm (\ref{zt})}, as $T\to\infty$, we have
\begin{equation}
\sup_{0\leq t\leq T}
\left|Z(t)-\frac12\left(\frac{J_1 \kappa_\alpha}{\sigma}
W_{1-\alpha/2}(t)\right)^2\right|=
o_{\AS}\left(T^{2-3\alpha/2+\alpha^2/4+\delta}+
T^{(1-\alpha/2+\gamma/2+\delta}\right),
\label{eq:main3}
\end{equation}
with $\gamma$ as in {\rm Theorem \ref{thm:main1}}, and arbitrary
$\delta>0$.
\end{thm}

\subsection{Consequences of Theorems \ref{thm:main1}--\ref{thm:main3}}
\label{sec:consequences}

First we deal with some immediate consequences of Theorem \ref{thm:main1}
for the partial sum process $S(t)$ as in (\ref{st}).

\begin{cor}\label{c1}
Under the conditions of {\rm Theorem \ref{thm:main1}}, as $n\to\infty$, we
have the following weak convergence in $D[0,\infty)$
(endowed with the uniform topology on compact sets).
\begin{equation}
\frac{\sigma(S(nt)-\mu nt)}{J_1 \kappa_\alpha n^{1-\alpha/2}}\convweak
W_{1-\alpha/2}(t).
\label{eq:weak1}
\end{equation}
\end{cor}

We note in passing that a more general first version of this
result was established directly on $D[0,1]$ by Taqqu \cite{Taqqu1975}.

\begin{cor}\label{c2}
Under the conditions of {\rm Theorem \ref{thm:main1}}, we have
\begin{equation}
\limsup_{n\to\infty}\frac{\sigma\sup_{0\leq t\leq 1}|S(nt)-\mu nt|}
{J_1 \kappa_\alpha n^{1-\alpha/2}(2\log\log n)^{1/2}}=1 \quad a.s.
\label{limsup}
\end{equation}
\begin{equation}
\liminf_{n\to\infty}\frac{\sigma\sup_{0\leq t\leq 1}|S(nt)-\mu nt|}
{J_1 \kappa_\alpha n^{1-\alpha/2}(\log\log n)^{-1+\alpha/2}}=c_\alpha
\quad a.s.,
\label{liminf}
\end{equation}
where $c_\alpha$ is a positive constant.
\end{cor}

\begin{cor}\label{c3}
Let $a_T$ be a nondecreasing function of $T$ such that
$T^{\tau+\delta}<a_T\leq T$ and $a_T/T$ is nonincreasing, where
$\delta>0$ is arbitrary, $\tau=(2-2\alpha)/(2-\alpha)$ if
$0<\alpha<1/2$ and $\tau=1/(2-\alpha)$ if $1/2\leq \alpha<1$.
Then, under the conditions of {\rm Theorem \ref{thm:main1}}, we have
\begin{equation}
\limsup_{T\to\infty}\frac{\sigma\sup_{0\leq t\leq T-a_T}\sup_{0\leq s\leq
a_T}|S(t+s)-S(t)-\mu s|}
{J_1 \kappa_\alpha a_T^{1-\alpha/2}(2(\log T/a_T+\log\log T))^{1/2}}=1
\quad a.s.
\label{limsupincr}
\end{equation}
\end{cor}

\begin{rem}\label{rem1}
{\rm As to the result of (\ref{limsup}), we note that, viewed via Theorem
\ref{thm:main1}, it is inherited from the LIL for fractional Brownian
motion that was established by Taqqu \cite{Taqqu1977}, cf.
Corollary A1, in a more general functional form. The latter can also be
spelled out for our $S(t)$ as in (\ref{st}), via Theorem \ref{thm:main1}.
Taqqu \cite{Taqqu1977} establishes his functional LIL for partial sums like
our (\ref{st}) and Gaussian processes as in his Theorem A1 separetely, on
their own. We note in passing that the statement of (21) also follows from
the LIL for Gaussian processes of Theorem 1.1 of Orey \cite{orey}.}
\end{rem}

\begin{rem}\label{rem2}
{\rm The result in (\ref{liminf}) is new for our partial sums, and it is
inherited from (3.6) of Theorem 3.3 of Monrad and Rootz\'en \cite{monrad},
where it is concluded for a standard fractional Brownian motion, with index
$1-\alpha/2$ in our terms. Thus the latter version of Chung's law of the
iterated logarithm for fractional Brownian motion is shared by partial
sums of LRD sequences of random variables. As to the constant $c_\alpha$,
its numerical value is not determined in \cite{monrad}.}
\end{rem}

\begin{rem}\label{rem3}
{\rm Ortega \cite{Ortega1984} extended the large increment results of
Cs\"org\H o--R\'ev\'esz \cite{CsR1979}, \cite{CsR1981} for a Wiener
process to centered Gaussian processes with stationary increments. Theorem
3 of \cite{Ortega1984} for a fractional Brownian motion as in
(\ref{covar}) reads as follows:}
For $T>0$ let $a_T$ be a nondecreasing function of $T$ such that
$0<a_T\leq T$ and $a_T/T$ is nonincreasing. Then
\begin{equation}
\limsup_{T\to\infty}\frac{\sup_{0\leq t\leq T-a_T}\sup_{0\leq s\leq a_T}
|W_{1-\alpha/2}(t+s)-W_{1-\alpha/2}(t)|}
{a_T^{1-\alpha/2}(2(\log T/a_T+\log\log T))^{1/2}}=1\quad a.s.
\label{ortega}
\end{equation}
{\rm Consequently, via Theorem \ref{thm:main1}, we arrive at
(\ref{limsupincr}). On taking $a_T=T$, we may conclude (\ref{limsup})
with $T$ instead of $n$.}
\end{rem}

We continue with spelling out new results for the counting process $N(t)$
as in (\ref{nt}) that follow from Theorem \ref{thm:main2}.

\begin{cor}\label{c4}
Under the conditions of {\rm Theorem \ref{thm:main2}}, as $n\to\infty$, we
have the following weak convergence in $D[0,\infty)$
(endowed with the uniform topology on compact sets).
\begin{equation}
\frac{\sigma(\mu N(\mu nt)-\mu nt)}{J_1 \kappa_\alpha
n^{1-\alpha/2}}\convweak
-W_{1-\alpha/2}(t).
\label{eq:weak2}
\end{equation}
\end{cor}

\begin{cor}\label{c5}
Under the conditions of {\rm Theorem \ref{thm:main2}}, we have
\begin{equation}
\limsup_{n\to\infty}\frac{\sigma\sup_{0\leq t\leq 1}|\mu N(\mu nt)-\mu
nt|}{J_1 \kappa_\alpha n^{1-\alpha/2}(2\log\log n)^{1/2}}=1 \quad a.s.
\label{limsup2}
\end{equation}
\begin{equation}
\liminf_{n\to\infty}\frac{\sigma\sup_{0\leq t\leq 1}|\mu N(\mu nt)-\mu
nt|}{J_1 \kappa_\alpha n^{1-\alpha/2}(\log\log n)^{-1+\alpha/2}}=c_\alpha
\quad a.s.,
\label{liminf2}
\end{equation}
where the positive constant $c_\alpha$ is that of (\ref{liminf}).
\end{cor}

\begin{cor}\label{c6}
Let $a_T$ be a nondecreasing function of $T$ such that
$T^{\tau+\delta}<a_T\leq T$ and $a_T/T$ is nonincreasing, where
$\delta>0$ is arbitrary, $\tau=(2-2\alpha)/(2-\alpha)$ if
$0<\alpha<1/2$ and $\tau=1/(2-\alpha)$ if $1/2\leq \alpha<1$. Then,
under the conditions of Theorem \ref{thm:main2}, we have
\begin{equation}
\limsup_{T\to\infty}\frac{\sigma\sup_{0\leq t\leq T-a_T}\sup_{0\leq s\leq
a_T}|\mu (N(\mu t+\mu s)-N(\mu t))-\mu s|}
{J_1 \kappa_\alpha a_T^{1-\alpha/2}(2(\log T/a_T+\log\log T))^{1/2}}=1
\quad a.s.
\label{limsupincr2}
\end{equation}
\end{cor}

\begin{rem}
{\rm We note that, mutatis mutandis, the conclusions of Corollaries
\ref{c4}, \ref{c5}, \ref{c6} for the counting process $N(t)$ follow via
Theorem \ref{thm:main2} exactly the same way as those of Corollaries
\ref{c1}, \ref{c2}, \ref{c3} do for $S(t)$ from Theorem \ref{thm:main1}
as noted in Remarks \ref{rem1}, \ref{rem2}, \ref{rem3}, i.e., from known
results for the fractional Brownian motion $W_{1-\alpha/2}(t)$.}
\end{rem}

The next corollaries deal with the process $Z(\cdot)$ as in (\ref{zt}) or,
equivalently, with the Vervaat process $V_n(\cdot)$ as in (\ref{Ve})  via
the strong approximation as in Theorem \ref{thm:main3}, in combination with
known results for the fractional Brownian motion $W_{1-\alpha/2}(t)$.

\begin{cor}\label{c7}
Under the conditions of {\rm Theorem \ref{thm:main3}}, as $n\to\infty$, we
have the following weak convergence in $C[0,\infty)$
(endowed with the uniform topology on compact sets).
\begin{equation}
\frac{2\sigma^2 Z(nt)}{J_1^2 \kappa_\alpha^2 n^{2-\alpha}}
=\frac{2\mu\sigma^2 V_n(t)}{J_1^2 \kappa_\alpha^2}\convweak
W_{1-\alpha/2}^2(t).
\label{eq:weak3}
\end{equation}
\end{cor}

\begin{cor}\label{c8}
Under the conditions of {\rm Theorem \ref{thm:main3}}, we have
\begin{equation}
\limsup_{n\to\infty}\frac{\sigma^2 \sup_{0\leq t\leq 1}Z(nt)}
{J_1^2 \kappa_\alpha^2 n^{2-\alpha}\log\log n}=
\limsup_{n\to\infty}\frac{\mu\sigma^2\sup_{0\leq t\leq 1}V_n(t)}
{J_1^2 \kappa_\alpha^2 \log\log n}=1 \quad a.s.
\label{limsup3}
\end{equation}
\begin{equation}
\liminf_{n\to\infty}\frac{2\sigma^2 \sup_{0\leq t\leq 1}Z(nt)}
{J_1^2 \kappa_\alpha^2 n^{2-\alpha}(\log\log n)^{-2+\alpha}}
=\liminf_{n\to\infty}\frac{2\mu \sigma^2 \sup_{0\leq t\leq 1}V_n(t)}
{J_1^2 \kappa_\alpha^2 (\log\log n)^{-2+\alpha}}=c_\alpha^2\quad a.s.,
\label{liminf3}
\end{equation}
where the positive constant $c_\alpha$ is as in (\ref{liminf}).
\end{cor}

\begin{rem}
{\rm If we were to assume the conditions of Theorem \ref{thm:main3} to
begin with, then, as noted already, we would have (\ref{eq:weak1})
directly via Taqqu \cite{Taqqu1975} that, in turn, in view of Theorem 1.1,
would lead to having (\ref{eq:weak2}) and (\ref{eq:weak3}) as well, as a
consequence of (\ref{eq:weak1}). Naturally, the
respective strong approximation results of Theorems \ref{thm:main1},
\ref{thm:main2}, \ref{thm:main3} were needed in order to conclude the
strong laws of Corollaries \ref{c2}, \ref{c5}, \ref{c8}, respectively.
We also note in passing that the strong approximation result of Theorem
\ref{thm:main2} leads to the weak convergence conclusion of Corollary
\ref{c4}, on assuming only condition (i) instead of condition (iii), that
is needed in the context of the first sentence above.}
\end{rem}

\begin{rem}
{\rm We note in passing that our positivity condition (iii) is used only
in Theorem \ref{thm:main3} in terms of our subordinated sequences
$Y_j=G(\tilde\eta_j)$, for which an arbitrary marginal distribution
function $F$ can always be obtained via choosing
$G(\cdot)=F^{-1}(\Phi(\cdot))$, where $\Phi(\cdot)$ stands for the standard
normal distribution function.}
\end{rem}

\begin{rem}
{\rm It would be of interest to prove analogues of our results at least
for functions $G$ whose Hermite rank is $m=2$. This would give rise to
the so-called Rosenblatt process (cf. Taqqu \cite{Taqqu1975}). The latter
process has stationary increments and covariance structure
like that of a fractional Brownian motion in (\ref{covar}). However, it is
non-Gaussian. Hence, difficulties arise when trying to deal with its
path behaviour, though one would think that it should be similar to that
of a fractional Brownian motion.}
\end{rem}

\section{Proofs}\label{sec:proofs}
To begin with, we note that the proof of Theorem \ref{thm:main3}, that is
based on Theorems \ref{thm:main1} and \ref{thm:main2}, also makes
fundamental use of the basic algebraic identity that is (2.5) of
\cite{CsakiCsorgoRychlikSteinebach2007}. In the i.i.d.
case, the appropriate approximations follow from
the Koml\'{o}s-Major-Tusn\'{a}dy results \cite{KomlosMajorTusnady1975},
\cite{KomlosMajorTusnady1976} and the work of Horv\'{a}th
\cite{Horvath1984}.

\subsection{Preliminary results and proof of Theorem
\ref{thm:main1}}\label{sec:proofmain1}

First we note that, via the proof of Lemma 6 in Oodaira
\cite{Oodaira1976}, we borrow the following construction
(cf. page 379 of Wang {\it et al.} \cite{WangLinGulati2003}) for our
fractional Brownian motion $W_{1-\alpha/2}$ that is being used throughout
this paper. We also note in passing that in the first line of page 379 in
\cite{WangLinGulati2003} one should write $\kappa_\alpha^2$ instead of
$\kappa_\alpha$.
\begin{lem}\label{lem1}
Let $\eta_j$ be defined by {\rm (\ref{linear})} with $\psi_k\sim
k^{-(1+\alpha)/2}$, $0<\alpha<1$, and put $\tilde\eta_j=\eta_j/\sigma$
with $\sigma^2:=\Exp(\eta_0^2)=\sum_{k=0}^\infty \psi_k^2$.
Then, on an appropriate probability space for the sequence
$\{\tilde\eta_j,\, j=0,1,\ldots\}$,
one can construct a fractional Brownian motion $W_{1-\alpha/2}(\cdot)$
such that, as $T\to\infty$, we have
\begin{equation}\label{eq:1}
\sup_{0\le t\le T}
\left|\kappa_{\alpha}^{-1}\sigma\sum_{j=1}^{[t]}\tilde\eta_j
-W_{1-\alpha/2}(t)\right|=o_{\AS}\left(T^{(1-\alpha)/2}\log T\right),
\end{equation}
where $\kappa_\alpha$ is as in {\rm (\ref{kalpha})}.
\end{lem}

We now state and prove a lemma, a strong reduction principle in terms of a
function with arbitrary Hermite rank. Earlier versions were given by Taqqu
\cite{Taqqu1977} and K\^ono \cite{kono1983}. The present version is of a
better rate that is based on combining a result of Lai and Stout
\cite{LaiStout1980} with that of Taqqu \cite{Taqqu1977}.

\begin{lem}\label{lem:CsSzW}
Assume that $H(\cdot)$ is an arbitrary function such that $\Exp
(H(\tilde\eta_0))=0$, $\Exp(H^2(\tilde\eta_0))<\infty$, and its Hermite
rank is $m\ge 1$, where $\tilde\eta_0$ is a standard normal random
variable. Let $\{\tilde\eta_j,\, j=0,1,\ldots\}$ be a stationary Gaussian
sequence with correlation as in $(\ref{cov})$. Then, as $n\to\infty$,
\begin{equation}
\sup_{0\le t\le n}\left|\sum_{j=1}^{[t]}H(\tilde\eta_j)-
\frac{J_m}{m!}\sum_{j=1}^{[t]}H_m(\eta_j)\right|=o_{\AS}(n^{\gamma/2
+\delta}),
\label{red}
\end{equation}
where $\gamma=2-(m+1)\alpha$ for $\alpha<1/(m+1)$, $\gamma=1$ for
$\alpha\ge 1/(m+1)$ and $\delta>0$ is arbitrary.
\end{lem}
{\bf Proof.}
Let
$$
U(n)=\sum_{j=1}^nH(\tilde\eta_j)-
\frac{J_m}{m!}\sum_{j=1}^nH_m(\tilde\eta_j).
$$
The Hermite rank of $H(\cdot)-J_m/m! H_m(\cdot)$ is at least $m+1$. By
Taqqu \cite[Proposition 4.2]{Taqqu1977} with $p=2$, we have for all
$a\geq 0$
$$
\Exp\left((U(n+a)-U(a))^2\right)\le
C\left\{n\sum_{i=0}^n|\rho_i|^{m+1}\right\}\le Cn^{\gamma+\delta},
$$
where $C$ is a finite positive constant, $\gamma$ is as above and
$\delta>0$ arbitrary. Consequently,

$$
P(|U(a+n)-U(a)|\geq x)\leq x^{-2}\Exp(U(a+n)-U(a))^2
\leq x^{-2}n^{\gamma+\delta}
$$
for $x>0$, $a\geq 0$, $n\geq 1$.

Since $\gamma+\delta>1$, the conditions of Theorem 7 of Lai and Stout
\cite{LaiStout1980} are satisfied with $g(n)=n^{\gamma+\delta}$, and
$p=2$. Hence, in our special case, with $\delta>0$ and $\gamma$ as
in Theorem \ref{thm:main1}, we conclude

$$
\lim_{n\to\infty}\frac{U(n)}{(n^{\gamma+\delta}(\log
n)^{1+\delta})^{1/2}}=0,\quad {\rm a.s.}
$$
Hence, as $n\to\infty$
$$
U(n)=o_{\AS}(n^{\gamma/2+\delta})
$$
that, in turn, also yields (\ref{red}) of Lemma 3.2.
\koniec

\medskip\noindent
{\bf Proof of Theorem \ref{thm:main1}}.
From Lemma \ref{lem:CsSzW} we obtain that, under the assumptions of
Theorem \ref{thm:main1}, as $T\to\infty$, we have the following rate in
the reduction principle for $\{Y_j=G(\tilde\eta_j);\, j=1,2,\ldots\}$:
\begin{equation}\label{eq:2}
\sup_{0\le t \le T}\left|\sum_{j=1}^{[t]}\left(G(\tilde\eta_j)-\mu\right)
-J_1\sum_{j=1}^{[t]}\tilde\eta_j\right|=o_{\AS}(T^{\gamma/2+\delta}),
\end{equation}
Since $1-\alpha< \gamma$, from (\ref{eq:1}) and (\ref{eq:2}) we
conclude that, as $T\to\infty$,
\begin{equation}\label{eq:3/2}
\sup_{0\le t\le T}\left|S(t)-\mu t-
\frac{J_1\kappa_{\alpha}}{\sigma}W_{1-\alpha/2}(t)\right|=
o_{\AS}(T^{\gamma/2+\delta}).
\end{equation}
\koniec

\subsection{Proof of Theorem \ref{thm:main2}}\label{sec:counting}
Assume the conditions of Theorem \ref{thm:main2}, i.e., those of Theorem
\ref{thm:main1} and condition (i).

We follow the approach of Horv\'ath \cite{Horvath1984},
\cite{Horvath1986}. We first show  that, as $T\to\infty$,
\begin{equation}
\sup_{0\leq t\leq T}|N(\mu t)-t|=o_{\AS}(T).
\label{nmut}
\end{equation}
Put $n=N(\mu t)$. Then $S(n)\leq \mu t\leq S(n+1)$. But according to
(\ref{eq:3/1}) we have
$$
|S(n)-\mu n|=o_{\AS}(n),\qquad |S(n+1)-\mu n|=o_{\AS}(n),
$$
from which
$$
|\mu t-\mu n|=o_{\AS}(n),
$$
i.e.
$$
t=N(\mu t)+o_{\AS}(N(\mu t)),
$$
which, in turn, implies (\ref{nmut}).

Let $\Lambda(t)$ and $\lambda(t)$ be two functions on $t\in [0,\infty)$
such that $\Lambda(u)=\inf\{t\ge 0:\lambda(t)>u\}$. Then (see, e.g.,
\cite{Horvath1984})
$$
\sup_{0<u<T}|\Lambda(u)-u|\le \sup_{0<t<\Lambda(T)}|\lambda(t)-t|.
$$
Consequently, as $T\to\infty$, (\ref{nmut}) in combination with
(\ref{limsup}) yields
$$
\sup_{0\leq t\leq T}|N(\mu t)-t| \leq
\sup_{0\leq t\leq N(\mu T)}|S(t)/\mu-t|
$$
\begin{equation}
\leq \sup_{0\leq t<(1+\varepsilon)T}|S(t)/\mu-t|
= o_{\AS}(T^{1-\alpha/2+\delta}).
\label{nmutmint}
\end{equation}

On writing now
\begin{equation}
\mu t-\mu N(\mu t)=(S(N(\mu t))-\mu N(\mu t))+(\mu t-S(N(\mu t))),
\label{mutminus}
\end{equation}
via (\ref{eq:3/1}) and (\ref{nmut}), for the first term of
(\ref{mutminus}), as
$T\to\infty$, we have uniformly in $t\in [0,T]$
\begin{eqnarray}\label{eq:6}
&& S(N(\mu t))-\mu N(\mu t)=\frac{J_1k_\alpha}{\sigma}
W_{1-\alpha/2}(N(\mu t))+o_{\AS}(T^{\gamma/2+\delta})\nonumber\\
&& = \frac{J_1 k_\alpha}{\sigma} W_{1-\alpha/2}(t)+
\frac{J_1 k_\alpha}{\sigma}(W_{1-\alpha/2}(N(\mu t))-W_{1-\alpha/2}(t))+
o_{\AS}(T^{\gamma/2+\delta}).~~~
\end{eqnarray}
Hence, it suffices to bound the increments
$$
\sup_{0\le t\le T}|W_{1-\alpha/2}(N(\mu t))-W_{1-\alpha/2}(t)|
$$
for large $T$. For doing this, we make use of a result of Ortega
\cite{Ortega1984} as quoted in (\ref{ortega}) of Remark 2.9.

On using now (\ref{ortega}) in combination with (\ref{nmut}), as
$T\to\infty$, we arrive at
\begin{eqnarray*}
&&\sup_{0\le t\le T}|W_{1-\alpha/2}(N(\mu t))-W_{1-\alpha/2}(t)|\\
&&
\leq \sup_{0\leq t\leq (1+\varepsilon)T} \sup_{0\leq s \leq
T^{1-\alpha/2}}
|W_{1-\alpha/2}(t+s)-W_{1-\alpha/2}(s)|
=O_{\AS}(T^{(1-\alpha/2)^2+\delta})
\end{eqnarray*}
with any $\delta>0$. Consequently, as $T\to\infty$, by (\ref{eq:6}), for
the first term of (\ref{mutminus}), we conclude
\begin{equation}\label{eq:11}
S(N(\mu t))-\mu N(\mu t)=
\kappa_{\alpha}J_1W_{1-\alpha/2}(t)+O_{\AS}(T^{(1-\alpha/2)^2
+\delta}+T^{\gamma/2+\delta})
\end{equation}
uniformly in $t\in [0,T]$.

As to the second term of (\ref{mutminus}), on putting $n=N(\mu t)$, we
have $S(n)\leq \mu t\leq S(n+1)$. Consequently, via (\ref{nmut}),
as $T\to\infty$, we have uniformly in $t\in [0,T]$
\begin{equation}
|S(N(\mu t))-\mu t|\leq \sup_{1\leq k\leq N(\mu t)}|S(k+1)-S(k)|
\leq \sup_{1\leq k\leq (1+\varepsilon)T}|S(k+1)-S(k)|
\label{snmu}
\end{equation}
with any $\varepsilon>0$. On estimating now the term on the right-hand
side of the last inequality in (\ref{snmu}) via (\ref{eq:3/1}), by
(\ref{snmu}) and using also (\ref{ortega}) with $a_T=1$, as $T\to\infty$,
we obtain
$$
\sup_{0\leq t\leq T}|S(N(\mu t))-\mu t)|
$$
\begin{equation}
\leq \frac{J_1 \kappa_\alpha}{\sigma}\sup_{0\leq s\leq
(1+\varepsilon)T}|W_{1-\alpha/2}(s+1)-W_{1-\alpha/2}(s)|
+o_{\AS}(T^{\gamma/2+\delta})=o_{\AS}(T^{\gamma/2+\delta}),
\label{secondterm}
\end{equation}
with any $\delta>0$.

On putting now together (\ref{secondterm}), (\ref{snmu}),
(\ref{mutminus}), the proof of Theorem \ref{thm:main2} is seen to be
complete.
\koniec

\subsection{Proof of Theorem \ref{thm:main3}}\label{sec:3-3}

Assume the conditions of Theorem \ref{thm:main3}, i.e., those of Theorem
\ref{thm:main1} and condition (iii). We note in passing that
condition (iii) implies that of (i) in view of our assumption throughout
that $J_1=\Exp(G(\tilde\eta_0)\tilde\eta_0)\neq 0$, for the sake of having the Hermite
rank of $G$ to be 1.

To begin with, we note that the algebraic identity of (2.5) in
\cite{CsakiCsorgoRychlikSteinebach2007} continues to hold true in our
present context. Consequently, with $S(t)$ and $Q(t)$ as in (\ref{st})
and (\ref{qt}) respectively, for $Z(t)$ as in (\ref{zt}), we have the
following identity:
\begin{equation}
Z(t)=\frac12 (S(t)-\mu t)^2+A(t)-\frac12 Q^2(t),
\label{zt1}
\end{equation}
where
$$
A(t)=\mu\int_{N(\mu t)}^t(S(s)-\mu s-(S(t)-\mu t))\, ds.
$$

In view of Theorem \ref{thm:main1} and (\ref{nmutmint}), when estimating
$A(t)$, we arrive at
\begin{equation}
A(t)=\frac{\mu \kappa_\alpha J_1}{\sigma}
\int_{N(\mu t)}^t(W_{1-\alpha/2}(s)-W_{1-\alpha/2}(t))\, ds
+o_{\AS}(T^{1-\alpha/2+\gamma/2+\delta}).
\label{at1}
\end{equation}
For estimating the integral in the latter conclusion, we have
$$
\left|\int_{N(\mu t)}^t(W_{1-\alpha/2}(s)-W_{1-\alpha/2}(t))\, ds\right|
$$
$$
\leq |N(\mu t)-t|\sup_{0\leq t\leq T}\sup_{0\leq u\leq
T^{1-\alpha/2+\delta}}|W_{1-\alpha/2}(t+u)-W_{1-\alpha/2}(t)|
$$
\begin{equation}
=o_{\AS}(T^{1-\alpha/2+(1-\alpha/2)^2+\delta}),
\label{at2}
\end{equation}
as $T\to\infty$, where we used (\ref{limsupincr2}) with $a_T=T$ in
Corollary \ref{c6}, as well as (\ref{ortega}) with $a_T=T^{1-\alpha/2}$
in Remark 2.9.
Thus, in view of (\ref{at1}) and (\ref{at2}), as $T\to\infty$, we arrive
at
\begin{equation}
\sup_{0\leq t\leq T}|A(t)|=o_{\AS}(T^{2-3\alpha/2+\alpha^2/4+\delta}
+T^{1-\alpha/2+\gamma/2+\delta}).
\label{at3}
\end{equation}

Next, in order to estimate $Q^2(t)$ in the identity (\ref{zt1}),
we have (cf. (\ref{qt}))
$$
Q(t)=S(t)-\mu t-(S(N(\mu t))-\mu N(\mu t))-(\mu t-S(N(\mu t))
$$
$$
=\frac{J_1\kappa_\alpha}{\sigma}W_{1-\alpha/2}(t)+o_{\AS}(t^{\gamma/2+\delta})
$$
$$
-\left(\frac{J_1\kappa_\alpha}{\sigma}W_{1-\alpha/2}(t)+
o_{\AS}(t^{\gamma/2+\delta})+o_{\AS}(t^{(1-\alpha/2)^2+\delta})\right)
$$
$$
+o_{\AS}(t^{\gamma/2+\delta}),
$$
where, as $t\to\infty$, we made respective use of (\ref{eq:3/1}) of
Theorem \ref{thm:main1}, (\ref{eq:11}) and (\ref{secondterm}).
Consequently, as $T\to\infty$, we arrive at
\begin{equation}
\sup_{0\leq t\leq T}Q^2(t)
=o_{\AS}(T^{\gamma+\delta}+T^{2(1-\alpha/2)^2+\delta}).
\label{qsquare}
\end{equation}

On combining now the identity of (\ref{zt1}) with (\ref{at3}) and
(\ref{qsquare}), we conclude Theorem \ref{thm:main3}.
\koniec

\bibliographystyle{plain}
\bibliography{lrd}

\begin{thebibliography}{10}

\bibitem{Bahadur1966}
R.~R. Bahadur.
\newblock A note on quantiles in large samples.
\newblock {\em Ann. Math. Statist.}, 37:577--580, 1966.

\bibitem{CsakiCsorgoKulik}
E. Cs{\'a}ki, M. Cs{\"o}rg{\H{o}} and R. Kulik.
\newblock On {V}ervaat processes for sums and renewals in weakly dependent 
cases.
\newblock In {\em {D}ependence in {P}robability, {A}nalysis and {N}umber
{T}heory}. A {V}olume in {M}emory of {W}alter {P}hilipp. Berkes et al., 
ed., pages 145--156. Kendrick Press, Heber City, UT, 2010.

\bibitem{CsakiCsorgoRychlikSteinebach2007}
E. Cs{\'a}ki, M. Cs{\"o}rg{\H{o}}, Z. Rychlik, and J. Steinebach.
\newblock On {V}ervaat and {V}ervaat-error-type processes for partial sums and
  renewals.
\newblock {\em J. Statist. Plann. Inference}, 137(3):953--966, 2007.

\bibitem{CsorgoKulik2008a}
M. Cs{\"o}rg{\H{o}} and R. Kulik.
\newblock Reduction principles for quantile and {B}ahadur-{K}iefer processes of
  long-range dependent linear sequences.
\newblock {\em Probab. Theory Related Fields}, 142(3-4):339--366, 2008.

\bibitem{CsorgoKulik2008b}
M. Cs{\"o}rg{\H{o}} and R. Kulik.
\newblock Weak convergence of {V}ervaat and {V}ervaat error processes of
  long-range dependent sequences.
\newblock {\em J. Theoret. Probab.}, 21(3):672--686, 2008.

\bibitem{CsR1979}
M. Cs{\"o}rg{\H{o}} and P. R\'ev\'esz.
\newblock How big are the increments of a Wiener process?
\newblock {\em Ann. Probab.} 7:731--737, 1979.

\bibitem{CsR1981}
M. Cs{\"o}rg{\H{o}} and P. R\'ev\'esz.
\newblock {\em Strong Approximations in Probability and Statistics.}
\newblock Academic Press, New York, 1981.

\bibitem{CsorgoSzyszkowiczWang2006}
M. Cs{\"o}rg{\H{o}}, B. Szyszkowicz, and L. Wang.
\newblock Strong invariance principles for sequential {B}ahadur-{K}iefer and
  {V}ervaat error processes of long-range dependent sequences.
\newblock {\em Ann. Statist.}, 34(2):1013--1044, 2006.

\bibitem{CSW2006cor}
M. Cs{\"o}rg{\H{o}}, B. Szyszkowicz, and L. Wang.
\newblock Correction: ``{S}trong invariance principles for sequential
  {B}ahadur-{K}iefer and {V}ervaat error processes of long-range dependent
  sequences'' [{A}nn. {S}tatist. {\bf 34} (2006), no. 2, 1013--1044;].
\newblock {\em Ann. Statist.}, 35(6):2815--2817, 2007.

\bibitem{DehlingTaqqu1989b}
H. Dehling and M. S. Taqqu.
\newblock The empirical process of some long-range dependent sequences with an
  application to {$U$}-statistics.
\newblock {\em Ann. Statist.}, 17(4):1767--1783, 1989.

\bibitem{Horvath1984}
L.~Horv{\'a}th.
\newblock Strong approximation of renewal processes.
\newblock {\em Stochastic Process. Appl.}, 18(1):127--138, 1984.

\bibitem{Horvath1986}
L.~Horv{\'a}th.
\newblock Strong approximations of renewal processes and their applications.
\newblock {\em Acta Math. Hungar.}, 47(1-2):13--28, 1986.

\bibitem{Kiefer1967}
J.~Kiefer.
\newblock On {B}ahadur's representation of sample quantiles.
\newblock {\em Ann. Math. Statist.}, 38:1323--1342, 1967.

\bibitem{Kiefer1970}
J.~Kiefer.
\newblock Deviations between the sample quantile process and the sample {${\rm
  df}$}.
\newblock In {\em Nonparametric {T}echniques in {S}tatistical {I}nference
  ({P}roc. {S}ympos., {I}ndiana {U}niv., {B}loomington, {I}nd., 1969)}, pages
  299--319. Cambridge Univ. Press, London, 1970.

\bibitem{KomlosMajorTusnady1975}
J.~Koml{\'o}s, P.~Major, and G.~Tusn{\'a}dy.
\newblock An approximation of partial sums of independent {${\rm RV}$}'s and
  the sample {${\rm DF}$}. {I}.
\newblock {\em Z. Wahrscheinlichkeitstheorie und Verw. Gebiete}, 32:111--131,
  1975.

\bibitem{KomlosMajorTusnady1976}
J.~Koml{\'o}s, P.~Major, and G.~Tusn{\'a}dy.
\newblock An approximation of partial sums of independent {RV}'s, and the
  sample {DF}. {II}.
\newblock {\em Z. Wahrscheinlichkeitstheorie und Verw. Gebiete}, 34(1):33--58,
  1976.

\bibitem{kono1983}
N. K{\^o}no.
\newblock Classical limit theorems for dependent random sequences having moment
  conditions.
\newblock In {\em Probability theory and mathematical statistics ({T}bilisi,
  1982)}, volume 1021 of {\em Lecture Notes in Math.}, pages 315--319.
  Springer, Berlin, 1983.

\bibitem{LaiStout1980}
T. L. Lai and W. Stout.
\newblock Limit theorems for sums of dependent random variables.
\newblock {\em Z. Wahrsch. Verw. Gebiete}, 51(1):1--14, 1980.

\bibitem{monrad}
D. Monrad and H. Rootz{\'e}n.
\newblock Small values of {G}aussian processes and functional laws of the
  iterated logarithm.
\newblock {\em Probab. Theory Related Fields}, 101(2):173--192, 1995.

\bibitem{oodaira1972}
H. Oodaira.
\newblock On {S}trassen's version of the law of the iterated logarithm for
  {G}aussian processes.
\newblock {\em Z. Wahrscheinlichkeitstheorie und Verw. Gebiete}, 21:289--299,
  1972.

\bibitem{Oodaira1976}
H. Oodaira.
\newblock Some limit theorems for the maximum of normalized sums of weakly 
dependent random variables. 
\newblock In Maruyama, G., and Prokhorov, J.V. (eds.) {\em Proceedings of 
the Third Japan-USSR Symposium on Probability Theory.} Lecture Notes in 
Math. Vol. 550, Springer, Berlin, Heidelberg, New York, pp. 467--474.

\bibitem{orey}
S. Orey.
\newblock Growth rate of certain {G}aussian processes.
\newblock In {\em Proceedings of the {S}ixth {B}erkeley {S}ymposium on
  {M}athematical {S}tatistics and {P}robability ({U}niv. {C}alifornia,
  {B}erkeley, {C}alif., 1970/1971), {V}ol. {II}: {P}robability theory}, pages
  443--451, Berkeley, Calif., 1972. Univ. California Press.

\bibitem{Ortega1984}
J. Ortega.
\newblock On the size of the increments of nonstationary {G}aussian processes.
\newblock {\em Stochastic Process. Appl.}, 18(1):47--56, 1984.

\bibitem{Taqqu1975}
M. S. Taqqu.
\newblock Weak convergence to fractional {B}rownian motion and to the
  {R}osenblatt process.
\newblock {\em Z. Wahrscheinlichkeitstheorie und Verw. Gebiete}, 31:287--302,
  1974/75.

\bibitem{Taqqu1977}
M. S. Taqqu.
\newblock Law of the iterated logarithm for sums of non-linear functions of
  {G}aussian variables that exhibit a long range dependence.
\newblock {\em Z. Wahrscheinlichkeitstheorie und Verw. Gebiete},
  40(3):203--238, 1977.

\bibitem{Taqqu2003}
M. S. Taqqu.
\newblock Fractional {B}rownian motion and long-range dependence.
\newblock In {\em Theory and applications of long-range dependence}, pages
  5--38. Birkh\"auser Boston, Boston, MA, 2003.

\bibitem{Vervaat1972b}
W.~Vervaat.
\newblock {\em Success epochs in {B}ernoulli trials (with applications in
  number theory)}.
\newblock Mathematisch Centrum, Amsterdam, 1972.
\newblock Mathematical Centre Tracts, No. 42.

\bibitem{Vervaat1972a}
W. Vervaat.
\newblock Functional central limit theorems for processes with positive drift
  and their inverses.
\newblock {\em Z. Wahrscheinlichkeitstheorie und Verw. Gebiete}, 23:245--253,
  1972.

\bibitem{WangLinGulati2003}
Q. Wang, Y-X. Lin, and C. M. Gulati.
\newblock Strong approximation for long memory processes with applications.
\newblock {\em J. Theoret. Probab.}, 16(2):377--389, 2003.

\end{thebibliography}

\end{document}